\documentclass{elsarticle}
\usepackage{amssymb, color}
\usepackage{amsmath}
 \usepackage{pifont}
 \usepackage{caption}
 \usepackage[utf8]{inputenc}
\usepackage{setspace}
\makeatletter
\def\LaTeX{\leavevmode L\raise.42ex
    \hbox{\kern-.3em\size{\sf@size}{0pt}\selectfont A}\kern-.15em\TeX}
\makeatother

\newcommand{\BibTeX}{{\rm B\kern-.05em{\sc
          i\kern-.025emb}\kern-.08em\TeX}}

\makeatletter
\def\@currentlabel{2.1}\label{e:dispaa}
\def\@currentlabel{2.21}\label{e:dispau}
\def\@currentlabel{2.22}\label{e:dispav}
\def\@currentlabel{2.23}\label{e:dispaw}
\def\@currentlabel{2.24}\label{e:dispax}
\def\theequation{\thesection.\@arabic\c@equation}
\makeatother

\renewcommand{\theequation}{\arabic{section}.\arabic{equation}}

\newcommand{\R}{\mathbb R}

\def \D{\Delta}

\newtheorem{thm}{Theorem} [section]
\newtheorem{lem}{Lemma} [section]

\newtheorem{cor}{Corollary} [section]

\newtheorem{rem}{Remark}[section]
\newtheorem{taggedtheoremx}{Theorem}
\newenvironment{taggedtheorem}[1]
 {\renewcommand\thetaggedtheoremx{#1}\taggedtheoremx}
 {\endtaggedtheoremx}

\renewcommand{\theequation}{\thesection.\arabic{equation}}
\renewcommand{\thesection}{\arabic{section}}
\renewcommand{\theequation}{\thesection.\arabic{equation}}
\let\ssection=\section\renewcommand{\section}{\setcounter{equation}{0}\ssection}

\begin{document}

\begin{frontmatter}

\title{On the classification of solutions to a weighted elliptic system
involving the Grushin operator.}
\author[fd]{Foued Mtiri}
\ead{mtirifoued@yahoo.fr}
\address[fd]{ANLIG, UR13ES32, University of Tunis El-Manar, 2092 El Manar II, Tunisia.}
\begin{abstract}
We investigate here the following weighted  degenerate elliptic system
\begin{align*}
-\Delta_{s} u =\Big(1+\|\mathbf{x}\|^{2(s+1)}\Big)^{\frac{\alpha}{2(s+1)}} v^p, \quad -\Delta_{s} v= \Big(1+\|\mathbf{x}\|^{2(s+1)}\Big)^{\frac{\alpha}{2(s+1)}}u^\theta,  \quad u,v>0\quad\mbox{in }\; \mathbb{R}^N:=\mathbb{R}^{N_1}\times \mathbb{R}^{N_2}.
\end{align*}
where  $\Delta_{s}=\Delta_{x}+|x|^{2s}\Delta_{y},$ is the Grushin operator, $s \geq 0,$  $\alpha \geq 0$ and $1<p\leq\theta.$ Here $$\|\mathbf{x}\|=\Big(|x|^{2(s+1)}+|y|^2\Big)^{\frac{1}{2(s+1)}}, \;\mbox{and}\;\; \mathbf{x}:=(x, y)\in \mathbb{R}^N:=\mathbb{R}^{N_1}\times \mathbb{R}^{N_2}.$$

In particular, we establish some new Liouville-type theorems for stable  solutions of the system,  which recover and considerably improve upon the known results \cite{cow, Hfh, HU, Fa, DP}. As a consequence, we obtain a
nonexistence result for the weighted Grushin equation
 \begin{align*}
 -\Delta_{s} u =\Big(1+\|\mathbf{x}\|^{2(s+1)}\Big)^{\frac{\alpha}{2(s+1)}} u^p,\;\; \quad u>0 \quad \mbox{in }\;\; \mathbb{R}^N.
\end{align*}

\end{abstract}
\end{frontmatter}
 \section{Introduction}
\setcounter{equation}{0}

We start by noting that throughout this article,  $N_s:=N_1+(1+s)N_2$ is called the
homogeneous dimension associated to the Grushin operator:
$$\Delta_{s}=\Delta_{x}+|x|^{2s}\Delta_{y},$$ where $s \geq 0,$ and
  $$\Delta_{x}:=\sum_{i=1}^{N_{1}}\frac{\partial^{2}}{\partial x_{i}^{2}},\quad\mbox{and} \;\; \Delta_{y}:=\sum_{j=1}^{N_{2}}\frac{\partial^{2}}{\partial y_{j}^{2}},$$ are Laplace operators with respect to $x\in \mathbb{R}^{N_1},$ $y\in \mathbb{R}^{N_2},$  and   $|x|^{2s}=\left(\sum \limits_{i=1}^{N_{1}}x_{i}^{2}\right)^{s}$.
\medskip

We devote this paper to the study of stable solutions to the following weighted  degenerate elliptic system

\begin{align}\label{1.1}
-\Delta_{s} u =\Big(1+\|\mathbf{x}\|^{2(s+1)}\Big)^{\frac{\alpha}{2(s+1)}} v^p, \quad -\Delta_{s} v= \Big(1+\|\mathbf{x}\|^{2(s+1)}\Big)^{\frac{\alpha}{2(s+1)}}u^\theta,  \quad u,v>0\quad\mbox{in }\; \mathbb{R}^N.
\end{align}
where $1<p\leq\theta,$ \; $\alpha \geq 0,$ and the  weighted Grushin equation
 \begin{align}\label{D}
 -\Delta_{s} u =\Big(1+\|\mathbf{x}\|^{2(s+1)}\Big)^{\frac{\alpha}{2(s+1)}} u^p, \quad u>0 \quad \mbox{in }\;\; \mathbb{R}^N ,
\end{align}
where $p>1,$ \; $\alpha \geq 0,$ and

$$\|\mathbf{x}\|=\Big(|x|^{2(s+1)}+|y|^2\Big)^{\frac{1}{2(s+1)}}, \quad s \geq 0, \quad\mbox{and}\quad \mathbf{x}:=(x, y)\in \mathbb{R}^N=\mathbb{R}^{N_1}\times \mathbb{R}^{N_2}, $$ is the norm corresponding to the Grushin distance, where $|x|$ and $|y|$ are the usual Euclidean norms in $\mathbb{R}^{N_1}$ and $\mathbb{R}^{N_2},$ respectively. Then, we can verify that the $\|\mathbf{x}\|$-norm is
1-homogeneous for the group of anisotropic dilations  attached to $\Delta_{s}.$ It is defined by
  $$\delta_\lambda(\mathbf{x})=(\lambda x, \lambda^{1+s}y), \quad \lambda>0 \quad\mbox{and}\quad \mathbf{x}:=(x, y)\in \mathbb{R}^N=\mathbb{R}^{N_1}\times \mathbb{R}^{N_2}.$$
  It is not hard to see that $ d \delta_\lambda(\mathbf{x})=\lambda^{N_1+(1+s)N_2}dxdy=\lambda^{N_s}d\mathbf{x},$ where $dxdy$ denotes the Lebesgue measure in $\mathbb{R}^N.$

\medskip
Our aim in this paper is to classify  stable solutions of  systems \eqref{1.1},  which can be regarded as a natural generalization 
 of the  weighted Grushin equation  \eqref{D}.  In order to state our results more accurately we defined  the notion of stability, we consider a general system given by
\begin{align}\label{1.gg2}
 -\Delta_{s} u = f(x,v),\quad -\Delta_{s} v= g(x,u),\quad \mbox{in }\; \mathbb{R}^N,
\end{align}
with $f,g\in C^1(\mathbb{R}^{N+1},\mathbb{R})$ satisfying $f_t:=\frac{\partial f(x,t)}{\partial t},
g_t=\frac{\partial f(x,t)}{\partial t}\geq 0$ in $\R.$ A smooth solution $(u,v)$ of \eqref{1.gg2} is said stable if there
exist positive smooth functions $\xi$, $\zeta$ verifying
\begin{align*}
 -\Delta_{s} \xi = f_{v}(x,v)\zeta, \quad -\Delta_{s} \zeta = g_{u}(x,u)\xi,\quad \mbox{in }\, \mathbb{R}^N.
\end{align*}
This definition is motivated by \cite{Mo, Fa, cow}. Let us review some results related to our problem.

  \medskip
Firstly, we recall the case  $s=\alpha=0$, the Lane–Emden equation and system have been extensively studied by many experts.
\smallskip

A celebrated result of Farina \cite{Far}, gives a complete classification up to the so called Joseph-Lundgren exponent (see \cite{CW, Far}).
 In  \cite{Fa}, using Farina’s approach, Fazly established Liouville type theorem of the weighted Lane-Emden equation \eqref{D} for $s = 0,$ $N$ satisfying $$N<N_{\bigstar}=2+\frac{2(2+\alpha)}{p-1}\left(p+\sqrt{p^2-p}\right),$$ and $p\geq 2.$ A large amount of works have been done
generalizing this result in various directions. To cite a few we refer to \cite{DY, DYfa, Hfh3, HU, Hfh}.

\medskip
We now turn to  the celebrated  Lane–Emden system  where $\alpha=s = 0,$  \eqref{1.1} reduces to
\begin{align}\label{1.2kk}
-\Delta u = v^p, \quad-\Delta v= u^\theta,\quad u,v>0\quad\mbox{in }\; \mathbb{R}^N, \quad\mbox{where }\;\theta \geq p > 1.
\end{align}
There is a famous conjecture who states that: {\sl Let $p, \theta > 0$. If the pair $(p, \theta)$ is subcritical, i.e.~if
\begin{align}\label{LjjE}
\frac{1}{p+1} + \frac{1}{\theta +1} > \frac{N-2}{N},
\end{align}
then there is no smooth solution to \eqref{1.2kk}.} This conjecture  is proved to be true for radial functions by Mitidieri \cite{em1}, (see also Serrin-Zou \cite{Zs1}), and for the nonradial solutions in dimensions $N = 3, 4$, by Souplet and his collaborators, see \cite{pqs, sou}. Other significant work in these topic can be found in \cite{s, ff, bm, Zs} and reference therein.

\medskip 

In  \cite{ MY}, the authors have
obtained the nonexistence of stable at infinity solutions of \eqref{1.2kk} for any $p, \theta > 0,$  satisfying \eqref{LjjE}.  We should also mention  that the classification for stable solutions of \eqref{1.2kk}, was obtained by Cowan \cite{cow}. Indeed, Using the stability inequality of system \eqref{1.2kk} and an interesting iteration argument, it was proved that there is no smooth stable solution to \eqref{1.2kk}, if $\max (1, 2t_0^-)<p\leq \theta,$   and $N$ satisfies $N<2+ \frac{4(\theta+1)}{p\theta-1}t_0^+,$ where \begin{align*}
t_{0}^{\pm} = \sqrt{\frac{p\theta(p+1)}{\theta+1}}\pm\sqrt{\frac{p\theta(p+1)}{\theta+1}-\sqrt{\frac{p\theta(p+1)}{\theta+1}}}.
\end{align*}
 In particular, if $N<10$, \eqref{1.2kk} has no stable solution for any $2<p\leq \theta.$ This result  was
 extended  in a   work  due to Hu \cite{ HU}, for the case of weighted  system  \eqref{1.1} with  $s=0,$
 namely
 \begin{align}\label{1.02kk}
-\Delta u = \rho(x) v^p, \quad-\Delta v= \rho(x) u^\theta,\quad u,v>0\quad\mbox{in }\; \mathbb{R}^N, \quad\mbox{where }\;\theta \geq p > 1,
\end{align}
and $\rho\equiv (1+|x|^2)^{\frac{\alpha}{2}},$  $\alpha>0.$  It was shown that, if  $2t_0^-<p\leq \theta$ and $N$ satisfies
$$N<2+ \frac{2(2+\alpha)(\theta+1)}{p\theta-1}t_0^+.$$
Furthermore, in  \cite{Hfh} the authors  established a Liouville type  theorem for the Lane-Emden system  with gereral
weights \eqref{1.02kk},  where  $\rho(x)$ is a radial function satisfying $\rho(x)\geq A(1+|x|^2)^{\frac{\alpha}{2}}$ at  infinity. This result was then
improved the previous works \cite{cow,Fa, HU}, and mainly obtained a new inverse comparison principle  which is the  key to deal with the case $1 < p \leq \frac{4}{3}.$  In particular, the range of nonexistence result in \cite{Hfh} is larger than that in \cite{cow, HU}.
\medskip

For the general  equation or  system with $s\neq 0$, the Liouville property is less understood and is more complicated to deal with than $s = 0.$  In the special  \textbf{case} $\alpha= 0,$ the system \eqref{1.1} and Eq \eqref{D}, become
\begin{align}\label{1b.1}
-\Delta_{s} u = v^p, \quad -\Delta_{s} v= u^\theta,  \quad u,v>0\quad\mbox{in }\; \mathbb{R}^N,
\end{align}
where $1<p\leq\theta,$  and  Grushin equation
 \begin{align}\label{Db}
 -\Delta_{s} u = u^p, \quad p>1,\quad u>0 \quad \mbox{in }\; \mathbb{R}^N.
\end{align}
\smallskip

Let us first recall some facts about the problem involving the Grushin operator. The  Liouville type theorem for  solutions of \eqref{Db}, has been established in  \cite{xxy, DDM} for  the case $1<p<\frac{N_s+2}{N_s-2}.$ The main tool in  \cite{xxy, DDM}, is the Kelvin transform combined with technique of moving planes.
  
  \smallskip

  Very recently, in  \cite{rahal} the author  extended some of Farina’s results  \cite{Far}  in order to prove  the nonexistence for nontrivial stable solution of the weighted  equation $-\Delta_{s} u=|x|_{s}^{\alpha} |u|^{p-1}u \;\;\; \mbox{in}\,\, \mathbb{R}^N,$ if

$$\left\{\begin{array}{ll}
1<p<\infty & \mbox{ si }\;  N_{s} \leq  10+4\alpha\\
1<p<p_{JL}(N_{s},\alpha)  & \mbox{ si }\; N_{s}>10+4\alpha,\end{array}
\right.$$
with
\begin{align*}
p_{JL}(N_{s},\alpha)=\frac{(N_{s}-2)^2 - 2(\alpha+2)(\alpha+N_{s}) +2\sqrt{(\alpha+2)^{3}(\alpha+2N_{s}-2)}}{(N_{s}-2)(N_{s}-10)}.
\end{align*}
 It should be notice that this condition $1<p<p_{JL}(N_{s},\alpha)$ is equivalent to $N_{s}<N_{\bigstar},$ where $N_{\bigstar}$ is given  explicitly in the above.

 \medskip
 Concerning  the nonexistence of  classical stable solutions of \eqref{1b.1} for  $\theta\geq p> 1,$  new results  were shown in  \cite{DP}, where the authors used the technics developed in  \cite{cow, Hfh} in order to obtain a direct extension of Theorem 1 in \cite{cow} for $s=0$:
 \begin{taggedtheorem}{A}
 \begin{enumerate}
\item Suppose that $\frac{4}{3}<p\leq \theta$ and
$$N_{s}<2+ \frac{4(\theta+1)}{p\theta-1}t_0^+,$$

Then there is no stable positive solution of \eqref{1b.1}. In particular, the assertion is true if $N_{s}\leq 10.$
\item Suppose that  $1< p \leq \min(\frac{4}{3}, \theta),$ and $$N_{s} < 2 + \left[2+\frac{2(p+1)}{p\theta-1}+\frac{4(2-p)}{\theta+p-2}\right]t_0^+.$$
    Then there is no bounded stable positive solution of  \eqref{1b.1}.
\end{enumerate}
\end{taggedtheorem}

\medskip
From now on, we assume that $\alpha > 0$. Our main objective is to prove  the following Liouville-type theorem  for 
stable solutions of \eqref{1.1} or \eqref{D} in $\R^N.$
\begin{thm}
\label{main3}
Let $x_0$ be the largest root of the polynomial
\begin{align}\label{newH}
H(x)=x^4 -
\frac{16p\theta(p+1)(\theta+1)}{(p\theta-1)^2}x^2 +
\frac{16p\theta(p+1)(\theta+1)(p+\theta+2)}{(p\theta-1)^3}x
-\frac{16p\theta(p+1)^2(\theta+1)^2}{(p\theta-1)^4}.
\end{align}
  \begin{enumerate}
\item
 If $\frac{4}{3}< p \leq \theta$ then \eqref{1.1} has no stable  solution if $N_{s}<2+(2+\alpha)x_0.$
 In particular, if $N_{s} \leq 10+4\alpha,$ then \eqref{1.1} has no  stable
solution for all $\frac{4}{3}<  p \leq \theta$.
\item If $1<p\leq \min(\frac{4}{3}, \theta)$, then \eqref{1.1} has no bounded  stable solution, if
$$N_{s} < 2 + \left[\frac{p}{2}+\frac{(2-p)(p \theta -1)}{(\theta+p-2)(\theta+1)}\right](\alpha+2)x_0.$$
\end{enumerate}
Therefore, if $N_{s} \leq 6+2\alpha,$ the system \eqref{1.1} has no bounded  stable
solution for all $\theta\geq p > 1.$
\end{thm}
\medskip

Thus, the following classification result for stable solution of \eqref{D}, is a direct consequence of Theorem \ref{main3}

\begin{cor}\label{main2}
\begin{enumerate}
 \item
 If $\frac{4}{3}< p$ then \eqref{D} has no stable  solution if
\begin{align}\label{N}
N_{s}<2+\frac{2(2+\alpha)}{p-1}\left(p+\sqrt{p^2-p}\right).
 \end{align}
 In particular, if $N_{s} \leq 10+4\alpha,$ then \eqref{D} has no stable 
solution for all $\frac{4}{3}< p$.
\item If $1<p\leq \frac{4}{3} $,
\eqref{D} has no bounded stable  solution for $N_{s}$ verifying \eqref{N}.
\end{enumerate}
Therefore, there is no bounded stable  solution of \eqref{D} for all $p > 1$ if $N_{s} \leq 10+4\alpha.$
\end{cor}
\medskip

\begin{rem} \begin{itemize}
\item If $s=0$, then the results in Theorem \ref{main3} and Corollary \ref{main2} coincide with that in \cite{Hfh}.
\item Using Remark 3.1 below, we see that  $2t_0^+\frac{\theta+1}{p\theta-1}<x_0$  for any  $1<p\leq \theta,$  where
$x_0$ is the largest root of the polynomial $H$ given by \eqref{newH}. So, Theorem \ref{main3} improves the bound given in Theorem {\bf A} with $\alpha=0.$
\end{itemize}
\end{rem}

 \medskip
 Establishing a Liouville type result for stable solution of \eqref{1.1}  is delicate, even we can borrow some ideas from  \cite{cow,Hfh, DP}.
 More precisely, the proof is based on nonlinear integral estimates  and the  comparison property between $u$ and $v.$ Nevertheless, unlike the  system \eqref{1.2kk}, the comparison property  has not been obtained  for the case of weighted system  \eqref{1.1}, since the operator $\Delta_{s}$ no longer has symmetry and it degenerates on the manifold $\{0\}\times \mathbb{R}^{N_2}$ and this introduces some essential difficulties in the proof of  Theorem \ref{main3}.  Furthermore,  to derive the  comparison property for the weighted Grushin operators,  we follow the general lines of the methods and techniques developed in  \cite{ DP}. Another difficulty, the $L^{1}$-estimates  to  the bootstrap  argument of   Cowan \cite{cow},   does not work in the case of  weighted Grushin operator. In order to overcome  the difficulties,  we instead pass to the $L^{2}$-estimates  in the bootstrap argument, which plays an essential role in iteration process.

\medskip
Our paper is organized as follows. In section 2, we prove some preliminaries results. The proofs of Theorem \ref{main3} and Corollary \ref{main2} are given in section 3.

 \section{Preliminary technical lemmas.}

\setcounter{equation}{0}
 In order to prove our results, we need some technical lemmas. In the following, $C$
denotes always a generic positive constant independent on $(u,v)$, which could be changed from one line to another.
The ball of center $0$ and radius $r > 0$ will be denoted by $B_r$.

\medskip
\subsection{Stability inequality}

\medskip
We can proceed similarly as the proof of Lemma 2.1 in \cite{DP}, we establish the following inequality

\begin{lem}\label{l.2.2}
If $(u,v)$ is a nonnegative stable solution of \eqref{1.1}, then
\begin{equation}
\label{1.3}  \sqrt{p\theta}\int_{\R^N} \Big(1+\|\mathbf{x}\|^{2(s+1)}\Big)^{\frac{\alpha}{2(s+1)}}
u^{\frac{\theta-1}{2}}v^{\frac{p-1}{2}}\phi^2dxdy  \leq
\int_{\R^N}|\nabla_{s}\phi|^2dxdy , \quad \forall\; \phi \in C_c^1(\R^N).
\end{equation}
\end{lem}

\noindent{\bf Proof.} Let $(u,v)$ denote a  stable solution of \eqref{1.1}. By the definition of stability, there exist positive smooth functions $\varphi$, $\psi$ verifying
\begin{align*}
 -\frac{\Delta_{s} \varphi}{ \varphi} = p \Big(1+\|\mathbf{x}\|^{2(s+1)}\Big)^{\frac{\alpha}{2(s+1)}}v^{p-1}\frac{\psi}{\varphi}, \quad -\frac{\Delta_{s} \psi}{ \psi} = \theta \Big(1+\|\mathbf{x}\|^{2(s+1)}\Big)^{\frac{\alpha}{2(s+1)}} u^{\theta-1}\frac{\varphi}{\psi}\quad \mbox{in }\, \R^N.
\end{align*}

Let $\gamma, \chi \in C_c^1(\R^N)$ and multiply the first equation by $\gamma^{2},$ and the second by $\chi^{2}$ and integrate
over $\R^N,$ to arrive at

\begin{align*}
 p \int_{\R^N}\Big(1+\|\mathbf{x}\|^{2(s+1)}\Big)^{\frac{\alpha}{2(s+1)}}v^{p-1}\frac{\psi}{\varphi}\gamma^{2}dxdy\leq
-\int_{\R^N}\frac{\Delta_{s} \varphi}{ \varphi}\gamma^{2}dxdy,
\end{align*}

 and

\begin{align*}
 \theta \int_{\R^N} \Big(1+\|\mathbf{x}\|^{2(s+1)}\Big)^{\frac{\alpha}{2(s+1)}} u^{\theta-1}\frac{\varphi}{\psi}\chi^{2}dxdy\leq
-\int_{\R^N}\frac{\Delta_{s} \psi}{ \psi}\chi^{2}dxdy.
\end{align*}

The simple calculation implies that
\begin{align*}
  \int_{\R^N}\left(-\frac{\Delta_{s} \varphi}{ \varphi}\gamma^{2}-|\nabla_{s}\gamma|^{2}\right)dxdy& \leq \int_{\mathbb{R}^N}\Big(\nabla_{s}\varphi\cdot \nabla_{s}(\gamma^{2}\varphi^{-1})-|\nabla_{s}\gamma|^{2}\Big)dxdy\\
& \leq \int_{\mathbb{R}^N}\Big(-\varphi^{-2}|\nabla_{s}\varphi|^{2}\gamma^{2}+2\varphi^{-1}\gamma\nabla_{s}\varphi\cdot\nabla_{s}\gamma-|\nabla_{s}\gamma|^{2}\Big)dxdy
\\
& \leq \int_{\mathbb{R}^N}-\Big(\varphi^{-1}\gamma\nabla_{s}\varphi-\nabla_{s}\gamma\Big)^{2}dxdy \leq 0
\end{align*}
 Proceeding as above,  we can easily show  that
\begin{align*}
-\int_{\R^N}\frac{\Delta_{s} \psi}{ \psi}\chi^{2}dxdy\leq\int_{\R^N} |\nabla_{s}\chi|^{2}dxdy.
\end{align*}

Using the inequality $2ab\leq a^{2}+b^{2},$   we deduce then

\begin{align*}
2\Big(1+\|\mathbf{x}\|^{2(s+1)}\Big)^{\frac{\alpha}{2(s+1)}} \sqrt{p\theta v^{p-1}{\varphi}u^{\theta-1}\gamma^{2}\chi^{2}}\leq \Big(1+\|\mathbf{x}\|^{2(s+1)}\Big)^{\frac{\alpha}{2(s+1)}} \Big( pv^{p-1}\frac{\psi}{\varphi}\gamma^{2}+\theta u^{\theta-1}\frac{\varphi}{\psi}\chi^{2}\Big).
\end{align*}

Taking $\varphi=\chi,$  and combining all these inequalities,  we get  readily the estimate \eqref{1.3}. \qed

\medskip
Inspired by the previous works \cite{sz, Fa, MPO}, we obtain the following integral
estimates for all  solutions of the system of \eqref{1.1}.

\begin{lem}\label{l.2.1}
Suppose that $(u,v)$ is a smooth solution of \eqref{1.1}, with $\theta \geq p > 1.$  Then
\begin{align}\label{VV}
  \int_{ B_{R}\times B_{R^{1+s}}} \Big(1+\|\mathbf{x}\|^{2(s+1)}\Big)^{\frac{\alpha}{2(s+1)}}v^{p} dxdy\leq C R^{N_{s}-\frac{2(\theta+1)p}{p\theta-1}-\frac{(p+1)\alpha}{p\theta-1}}.
\end{align}
\begin{align}\label{UU}
   \int_{ B_{R}\times B_{R^{1+s}}}\Big(1+\|\mathbf{x}\|^{2(s+1)}\Big)^{\frac{\alpha}{2(s+1)}} u^{\theta} dxdy \leq C R^{N_{s}-\frac{2(p+1)\theta}{p\theta-1}-\frac{(\theta+1)\alpha}{p\theta-1}}.
\end{align}
\end{lem}

\noindent{\bf Proof.} Let $\chi_{j}\in C_c^\infty\left(\mathbb{R}, [0,1]\right),$ $j=1,2$ be a cut-off function verifying $0 \leq \chi_{j} \leq 1,$ $$\chi_{j}=1 \quad \mbox{on} [-1,1],\quad \mbox{and} \quad  \chi_{j}=0 \quad \mbox{outside } \quad [-2^{1+(j-1)s},2^{1+(j-1)s}].$$
For $R\geq 1,$ put $\psi_{R}(x,y)=\chi_{1}(\frac{x}{R})\chi_{2}(\frac{y}{R^{1+s}}),$ it is easy to verify  that there exists $C >0$ independent of $R$ such that

$$|\nabla_{x}\psi_{R}|\leq \frac{C}{R}\quad \mbox{and} \quad  |\nabla_{y}\psi_{R}|\leq \frac{C}{R^{1+s}},$$

\smallskip

$$|\D_{x}\psi_{R}|\leq \frac{C}{R^{2}}\quad \mbox{and} \quad  |\D_{y}\psi_{R}|\leq \frac{C}{R^{2(1+s)}}.$$

 Multiplying the equation $-\D_{s} u=\Big(1+\|\mathbf{x}\|^{2}\Big)^{\frac{\alpha}{2}}v^{p}$ by $\psi ^{m},$ and integrating by parts to yield
\begin{align*}
& \;\int_{ B_{2R}\times B_{(2R)^{1+s}}} \Big(1+\|\mathbf{x}\|^{2(s+1)}\Big)^{\frac{\alpha}{2(s+1)}}v^{p}\psi_{R}^{m} dxdy\\
& \;= -\int_{ B_{2R}\times B_{(2R)^{1+s}}}u\Delta_{s}(\psi_{R}^{m}) dxdy \leq \frac{C}{R^{2}}\int_{ B_{2R}\times B_{(2R)^{1+s}}} u \psi_{R}^{m-2} dxdy.
\end{align*}
 Let $\frac{1}{\theta}+\frac{1}{\theta'}=1.$  Apply H\"older's inequality,  we obtain
\begin{align*}
& \; \int_{ B_{2R}\times B_{(2R)^{1+s}}} \Big(1+\|\mathbf{x}\|^{2(s+1)}\Big)^{\frac{\alpha}{2(s+1)}}v^{p}\psi_{R}^{m} dxdy\\
\leq&\;\frac{C}{R^{2}}\left(\int_{ B_{2R}\times B_{(2R)^{1+s}}}\Big(1+\|\mathbf{x}\|^{2(s+1)}\Big)^{-\frac{\theta'\alpha}{2(s+1)\theta}} dxdy\right)^{\frac{1}{\theta'}}\\
&\; \;\;\;\times\left(\int_{ B_{2R}\times B_{(2R)^{1+s}}}\Big(1+\|\mathbf{x}\|^{2(s+1)}\Big)^{\frac{\alpha}{2(s+1)}} u^{\theta} \psi_{R}^{(m-2)\theta}dxdy\right)^{\frac{1}{\theta}}\\
\leq& \; C R^{\frac{N_{s}}{\theta'} -\frac{\alpha}{\theta}-2}\left(\int_{ B_{2R}\times B_{(2R)^{1+s}}}\Big(1+\|\mathbf{x}\|^{2(s+1)}\Big)^{\frac{\alpha}{2(s+1)}} u^{\theta} \psi_{R}^{(m-2)\theta}dxdy\right)^{\frac{1}{\theta}}.
\end{align*}
Adopting the similar argument as  above where we use  the second equation in \eqref{1.1}, we deduce then  for $k \geq 2$

\begin{align*}
  & \; \int_{ B_{2R}\times B_{(2R)^{1+s}}} \Big(1+\|\mathbf{x}\|^{2(s+1)}\Big)^{\frac{\alpha}{2(s+1)}}u^{\theta}\psi_{R}^{k} dxdy \\
   \leq& \;C R^{\frac{N_{s}}{p'} -\frac{\alpha}{p}-2}\left(\int_{ B_{2R}\times B_{(2R)^{1+s}}}\Big(1+\|\mathbf{x}\|^{2(s+1)}\Big)^{\frac{\alpha}{2(s+1)}} v^{p} \psi_{R}^{(k-2)p}dxdy\right)^{\frac{1}{p}},
\end{align*}
where $\frac{1}{p}+\frac{1}{p'}=1.$ Take now $k$ and $m$ large verifying $ m \leq (k-2)p$ and $k \leq (m-2)\theta.$ Combining the two
above inequalities, one concludes
\begin{align*}
  & \;\int_{ B_{2R}\times B_{(2R)^{1+s}}} \Big(1+\|\mathbf{x}\|^{2(s+1)}\Big)^{\frac{\alpha}{2(s+1)}}v^{p}\psi_{R}^{m} dxdy\\
   \leq & \; C R^{\frac{N_{s}}{\theta'} -\frac{\alpha}{\theta}-2} R^{\big(\frac{N_{s}}{p'} -\frac{\alpha}{p}-2\big)\frac{1}{\theta}}\left(\int_{ B_{2R}\times B_{(2R)^{1+s}}}\Big(1+\|\mathbf{x}\|^{2(s+1)}\Big)^{\frac{\alpha}{2(s+1)}} v^{p} \psi_{R}^{(k-2)p}dxdy\right)^{\frac{1}{p\theta}}\\
   \leq & \; C R^{N_{s} - \frac{N_{s}}{p\theta} -\frac{\alpha(p+1)}{p\theta}- \frac{2(\theta+1)}{\theta}}\left(\int_{ B_{2R}\times B_{(2R)^{1+s}}}\Big(1+\|\mathbf{x}\|^{2(s+1)}\Big)^{\frac{\alpha}{2(s+1)}} v^{p} \psi_{R}^m dxdy\right)^{\frac{1}{p\theta}}.
\end{align*}
So, we get
\begin{align*}
 & \; \int_{ B_{R}\times B_{R^{1+s}}} \Big(1+\|\mathbf{x}\|^{2(s+1)}\Big)^{\frac{\alpha}{2(s+1)}}v^{p}dxdy\\
   \leq & \; \int_{ B_{2R}\times B_{(2R)^{1+s}}} \Big(1+\|\mathbf{x}\|^{2(s+1)}\Big)^{\frac{\alpha}{2(s+1)}}v^{p}\psi_{R}^{m} dxdy
   \leq  \; C R^{N_{s}-\frac{2(\theta+1)p}{p\theta-1}-\frac{(p+1)\alpha}{p\theta-1}}.
\end{align*}
Finally,  the estimate \eqref{UU} follows from using the same argument as above.   \qed

\medskip
\subsection{Comparison principle}

\medskip
Inspired by the previous works in \cite{DP,Hfh}, we can find the point-wise estimate of solution (u, v).

\begin{lem}\label{Soup}(Comparaison property.) Let $\theta \geq p > 1,$ $\alpha\geq 0,$ and  $(u, v)$ be 
positive solution of  \eqref{1.1}. Then there holds
\begin{align}
\label{estS}
u^{\theta + 1} \leq \frac{\theta+1}{p+1}v^{p+1} \;\; \mbox{in}\;\;\mathbb{R}^N.
\end{align}
\end{lem}
\noindent
{\bf Proof.} Let  $\sigma=\frac{p+1}{\theta+1}\in (0,1],$ $\lambda = \sigma^\frac{-1}{\theta+1},$ and $w= u- \lambda v^{\sigma}.$ A straightforward computation implies
\begin{align*}
\D_{s} w = \D_{s} u-\lambda\sigma v^{\sigma-1} \D_{s}v -\lambda\sigma (\sigma-1)|\nabla_{s}v|^{2}v^{\sigma-2}&\geq\Big(1+\|\mathbf{x}\|^{2(s+1)}\Big)^{\frac{\alpha}{2(s+1)}}\left[ - v^p+\lambda\sigma v^{\sigma-1}u^{\theta} \right]\\
&=\Big(1+\|\mathbf{x}\|^{2(s+1)}\Big)^{\frac{\alpha}{2(s+1)}}v^{\sigma-1}\left[ - v^{p-\sigma+1}+\lambda\sigma u^{\theta} \right]\\
&=\Big(1+\|\mathbf{x}\|^{2(s+1)}\Big)^{\frac{\alpha}{2(s+1)}}v^{\sigma-1}\left[ \lambda^{-\theta}u^{\theta} - v^{\theta\sigma} \right].
\end{align*}

Therefore, for any $\sigma\in (0,1]$, there exists $C > 0$ such that
\begin{align}\label{xy}
Cv^{\sigma-1}\left[u^{\theta} - \Big(\lambda v^{\sigma} \Big)^{\theta}\right]\leq\Big(1+\|\mathbf{x}\|^{2(s+1)}\Big)^{\frac{\alpha}{2(s+1)}}v^{\sigma-1}\left[\frac{u^{\theta} - \Big(\lambda v^{\sigma} \Big)^{\theta}}{\lambda^{\theta}} \right]\leq \D_{s} w.
\end{align}

We need to prove that $$u\leq\lambda v^{\sigma}.$$
 We shall show that $w \leq0,$  by a contradiction argument. Suppose that
 \begin{align}\label{00xy}
\sup_{\mathbb{R}^N}w>0.
\end{align}

 Next, we split the proof into two cases.

\medskip
\textit{Case 1: } we consider the case where the supremum of $w$ is attained at infinity.

\medskip
Choose now   $\phi_{R}(x,y)=\psi^{m}(\frac{x}{R},\frac{y}{R^{1+s}}),$  where $m>0,$ and $\psi $ a cut-off function in $ C_c^\infty\left(\mathbb{R}^N, [o,1]\right),$ such that $$\psi=1 \quad \mbox{on} \quad B_{1}\times B_{1},\quad \mbox{and} \quad  \psi=0 \quad \mbox{outside } \quad B_{2}\times B_{2^{1+s}}.$$

A simple calculation,  implies that

$$\frac{|\nabla_{s}\phi_{R}|^{2}}{\phi_{R}}\leq  \frac{C}{R^{2}}\phi_{R}^{\frac{m-2}{2}}\quad \mbox{and} \quad  |\D_{s}(\phi_{R})|\leq \frac{C}{R^{2}}\phi_{R}^{\frac{m-2}{2}}.$$

Set $$w_{R}=\phi_{R}w,$$ which is a compactly supported function. Then there exists $(x_{R},y_{R})\in  B_{2R}\times B_{(2R)^{1+s}},$ such that
\begin{align*}
 w_{R}(x_{R},y_{R})=\max_{\mathbb{R}^N}w_{R}(x,y)\rightarrow \sup_{\mathbb{R}^N}w(x,y) \quad \mbox{as }\; R\rightarrow \infty,
\end{align*}
which implies
\begin{align*}
\nabla_{s} w_{R}(x_{R},y_{R})=0 \quad \mbox{and}\quad \D_{s} w_{R}(x_{R},y_{R})\leq0,
\end{align*}
which means that at $(x_{R},y_{R})$,
\begin{align}\label{2xy}
\nabla_{s} w=-\phi_{R}^{-1}\nabla_{s}\phi_{R}w \quad \mbox{and}\quad \phi_{R}\D_{s} w\leq 2w\phi_{R}^{-1}|\nabla_{s}\phi_{R}|^{2}-w\D_{s}\phi_{R}.
\end{align}
 From \eqref{2xy},  and using the properties of $\phi_{R}$,  we can conclude then

 \begin{align}\label{3xy}
\phi_{R}\D_{s} w\leq \frac{C}{R^{2}}\phi_{R}^{\frac{m-2}{2}}w.
\end{align}

 Furthermore, for  $w= u- \lambda v^{\sigma}\geq0,$ we  observe that
 \begin{align}\label{4xy}
 \frac{u^{\theta}}{w^{\theta}}-\frac{(\lambda v^{\sigma})^{\theta}}{w^{\theta}}\geq1, \quad \mbox{ or equivalently }\quad
\lambda^{-\theta}u^{\theta} - v^{\theta\sigma}\geq\lambda^{-\theta} w^{\theta}.
\end{align}
 Multiplying  \eqref{xy} by $\phi_{R},$ combining  it with \eqref{4xy} and \eqref{3xy}, one obtains

  \begin{align*}
v^{\sigma-1} w_{R}^{\theta}\phi_{R}^{\frac{m+2}{2}}\leq \frac{C}{R^{2}}w_{R}\phi_{R}.
\end{align*}
As $\sigma\leq1.$ If the sequence $v(x_{R},y_{R})$ is bounded, and we choose
 \begin{align*}
\theta=\frac{m+2}{m} \quad \mbox{so that }\;m=\frac{2}{\theta-1},
\end{align*}
there holds then
   \begin{align*}
 w_{R}^{\theta-1}(x_{R},y_{R})\leq \frac{C}{R^{2}}.
\end{align*}
 Taking the limit $R \to\infty$, we have $\sup_{\mathbb{R}^N}w=0,$ which contradicts \eqref{00xy}, the claim follows.

\medskip
\textit{Case 2: } If there exists $(x^{0}, y^{0}),$ such that $\sup_{\mathbb{R}^N}w(x^{0}, y^{0})=u(x^{0}, y^{0})- \lambda v^{\sigma}(x^{0}, y^{0})>0,$ then

$\frac{\partial w}{\partial x}(x^{0}, y^{0})=0,$  or  $\frac{\partial w}{\partial y}(x^{0}, y^{0})=0,$ and  $\frac{\partial^{2} w}{\partial x^{2}}(x^{0}, y^{0})\leq0,$  or  $\frac{\partial^{2} w}{\partial y^{2}}(x^{0}, y^{0})\leq0.$  This contradicts the fact that  $\D_{s} w(x^{0}, y^{0})\geq0,$  the proof is completed. \qed

\medskip
Exploiting the technique introduced in  \cite{Hfh}, we establish  the following comparison property between $u$ and $v$, which is
somehow an inverse version of the  estimate  \eqref{estS}.

\begin{lem}\label{invSoup}
Suppose that $(u, v)$ be a  solution of (1.1), with $\theta \geq p > 1$ and assume that $u$ is bounded, then
 \begin{align}\label{investS}
v\leq \|u\|_{\infty}^\frac{\theta-p}{p+1}u.
 \end{align}
\end{lem}
\noindent
{\bf Proof.} Let $w= v- \l u,\;$ where $\l = \|u\|_{\infty}^\frac{\theta-p}{p+1}$. We have, as $\theta \geq p$,
\begin{align}\label{x1y}
 \begin{split}
\D w = \Big(1+\|\mathbf{x}\|^{2(s+1)}\Big)^{\frac{\alpha}{2(s+1)}}\left(\l v^p- u^\theta\right) &=\Big(1+\|\mathbf{x}\|^{2}\Big)^{\frac{\alpha}{2}}\left[ \l v^p- \left(\frac{u}{ \|u\|_{\infty}}\right)^{\theta} \|u\|_{\infty}^{\theta}\right] ]\\&\geq \Big(1+\|\mathbf{x}\|^{2(s+1)}\Big)^{\frac{\alpha}{2(s+1)}}\left[\l v^p- \left(\frac{u}{ \|u\|_{\infty}}\right)^p \|u\|_{\infty}^{\theta}\right]\\
&=\Big(1+\|\mathbf{x}\|^{2(s+1)}\Big)^{\frac{\alpha}{2(s+1)}} \|u\|_{\infty}^{\theta - p} \left(\frac{\l v^p}{\|u\|_{\infty}^{\theta - p}} - u^p\right)\\
&= \Big(1+\|\mathbf{x}\|^{2(s+1)}\Big)^{\frac{\alpha}{2(s+1)}}\|u\|_{\infty}^{\theta - p} \left(\l^{-p}v^p - u^p\right).
\end{split}
\end{align}
For the rest of the proof, can be proceeded  as in Lemma 2.3 where we replace just  \eqref{xy} by  \eqref{x1y}, so we omit the details.  \qed

\medskip
A crucial ingredient in our proof of  Theorem \ref{main3} for   the case $1 < p \leq \frac{4}{3},$ is given by the following Lemma is a consequence of the stability inequality \eqref{1.3},  and the inverse comparison principle \eqref{investS}.
\begin{lem}
\label{lemnew}
Let $(u,v)$ be a stable solution to \eqref{1.1} with $1< p \leq \min(\frac{4}{3}, \theta)$. Assume that $u$ is bounded, there holds
\begin{align}
\label{2.3}
\int_{ B_{R}\times B_{R^{1+s}}}\Big(1+\|\mathbf{x}\|^{2(s+1)}\Big)^{\frac{\alpha}{2(s+1)}} v^2dxdy \leq  CR^{N_{s}-\frac{2(\theta+1)p}{p\theta-1}-\frac{(p+1)\alpha}{p\theta-1} - \frac{2(2+\alpha)(2 - p)}{\theta+p-2}}.
\end{align}
\end{lem}

\noindent{\bf Proof.} Let $(u,v)$ be a stable solution of \eqref{1.1}, where $u$ is bounded. Take $\chi \in  C_c^\infty\left(\mathbb{R}^N, [0,1]\right),$ be a cut-off function satisfying  $\chi=1 \; \mbox{on} \; B_{1}\times B_{1},\; \mbox{and} \;  \chi=0 \; \mbox{outside } \; B_{2}\times B_{2^{1+s}}.$  Put $\eta_{R}(x,y)=\chi(\frac{x}{R},\frac{y}{R^{1+s}}).$

 Thanks to the approximation argument, the stability property  \eqref{1.3} holds true with $\phi = v \eta_{R}$, we deduce then
\begin{align}\label{x7y}
 \begin{split}
 & \;\sqrt{p\theta}\int_{\mathbb{R}^N} \Big(1+\|\mathbf{x}\|^{2(s+1)}\Big)^{\frac{\alpha}{2(s+1)}}u^{\frac{\theta-1}{2}} v^{\frac{p-1}{2}}v^2 \eta_{R}^2dxdy\\
& \leq  \; \int_{\mathbb{R}^N} |\nabla_{s} v|^2 \eta_{R}^2dxdy +  \int_{\mathbb{R}^N}v^2 |\nabla_{s} \eta_{R} |^2dxdy -\frac{1}{2}\int_{\mathbb{R}^N}v^2 \D_{s} (\eta_{R}^2)dxdy.
 \end{split}
\end{align}

Multiplying $-\D v= \Big(1+\|\mathbf{x}\|^{2(s+1)}\Big)^{\frac{\alpha}{2(s+1)}} u^{\theta}$ by $v \eta_{R} ^2$ and integrating by parts,  one gets:
\begin{align*}
  \int_{\mathbb{R}^N} |\nabla_{s} v|^2 \eta_{R}^2dxdy= \int_{\mathbb{R}^N}\Big(1+\|\mathbf{x}\|^{2(s+1)}\Big)^{\frac{\alpha}{2(s+1)}} u^\theta v \eta_{R}^2dxdy + \frac{1}{2}\int_{\mathbb{R}^N}v^2 \D_{s} (\eta_{R}^2)dxdy.
\end{align*}
 From \eqref{estS}, ones has
\begin{align*}
  \int_{\mathbb{R}^N} |\nabla_{s} v|^2 \eta_{R}^2dxdy\leq  \sqrt{\frac{\theta+1}{p+1}}\int_{\mathbb{R}^N}\Big(1+\|\mathbf{x}\|^{2(s+1)}\Big)^{\frac{\alpha}{2(s+1)}} u^{\frac{\theta-1}{2}} v^{\frac{p+1}{2}}v \eta_{R}^2dxdy + \frac{1}{2}\int_{\mathbb{R}^N}v^2 \D_{s} (\eta_{R}^2)dxdy.
\end{align*}
Substituting this in \eqref{x7y}, we obtain readily

\begin{align*}
 \left(\sqrt{p \theta}-\sqrt{\frac{\theta+1}{p+1}}\right)\int_{\mathbb{R}^N}\Big(1+\|\mathbf{x}\|^{2(s+1)}\Big)^{\frac{\alpha}{2(s+1)}} u^{\frac{\theta-1}{2}} v^{\frac{p+3}{2}} \eta_{R}^2dxdy \leq   \int_{\mathbb{R}^N}v^2 |\nabla_{s} \eta_{R} |^2dxdy,
\end{align*}

 Take $\eta_{R} = \varphi_{R}^{m}$ with $m > 2.$ Using Lemma \ref{invSoup}, there exists a positive constant $C$ such that

 \begin{align}\label{x8y}
 \begin{split}
 \int_{\mathbb{R}^N}\Big(1+\|\mathbf{x}\|^{2(s+1)}\Big)^{\frac{\alpha}{2(s+1)}}v^{\frac{\theta+p+2}{2}}\varphi_{R}^{2m}dxdy \leq \frac{C}{R^{2+\alpha}}\int_{\mathbb{R}^N}\Big(1+\|\mathbf{x}\|^{2(s+1)}\Big)^{\frac{\alpha}{2(s+1)}} v^2\varphi_{R}^{2m-2}dxdy.
 \end{split}
 \end{align}
Since $1< p \leq \min(\frac{4}{3}, \theta),$ a direct calculation shows that

\begin{align*}
2 = p\lambda + \frac{\theta + p + 2}{2}(1 - \lambda) \quad \mbox{with } \lambda = \frac{\theta + p - 2}{\theta + 2 - p} \in (0, 1).
 \end{align*}
By H\"older's inequality, the integral  in the right hand side of \eqref{x8y}, can be estimated as
  \begin{align}\label{x9y}
 \begin{split}
&\int_{\mathbb{R}^N}\Big(1+\|\mathbf{x}\|^{2(s+1)}\Big)^{\frac{\alpha}{2(s+1)}} v^2\varphi_{R}^{2m-2}dxdy\\
& \leq\left( \int_{\mathbb{R}^N}\Big(1+\|\mathbf{x}\|^{2(s+1)}\Big)^{\frac{\alpha}{2(s+1)}}v^{\frac{\theta+p+2}{2}}\varphi_{R}^{2m}dxdy\right)^{1 - \lambda} \left(\int_{\R^N}\Big(1+\|\mathbf{x}\|^{2(s+1)}\Big)^{\frac{\alpha}{2(s+1)}}v^p\varphi_{R}^{2m\lambda - 2}dxdy\right)^{\lambda}\\
&\leq\left( \int_{\mathbb{R}^N}\Big(1+\|\mathbf{x}\|^{2(s+1)}\Big)^{\frac{\alpha}{2(s+1)}}v^{\frac{\theta+p+2}{2}}\varphi_{R}^{2m}dxdy\right)^{1 - \lambda} \left(\int_{ B_{2R}\times B_{(2R)^{1+s}}}\Big(1+\|\mathbf{x}\|^{2(s+1)}\Big)^{\frac{\alpha}{2(s+1)}}v^pdxdy\right)^{\lambda}.
\end{split}
 \end{align}

 For the last inequality, we used $0\leq\varphi_{R}\leq1$ and chosen $m$ large such that $m\lambda > 1.$ Combining \eqref{VV} and  \eqref{x8y}-\eqref{x9y}, there holds
  \begin{align*}
&\;\int_{\mathbb{R}^N}\Big(1+\|\mathbf{x}\|^{2(s+1)}\Big)^{\frac{\alpha}{2(s+1)}} v^2\varphi_{R}^{2m-2}dxdy \\
& \leq  \; R^{\lambda\left(N_{s}-\frac{2(\theta+1)p}{p\theta-1}-\frac{(p+1)\alpha}{p\theta-1}\right)-(2+\alpha)(1 - \lambda)}\left( \int_{\mathbb{R}^N}\Big(1+\|\mathbf{x}\|^{2(s+1)}\Big)^{\frac{\alpha}{2(s+1)}} v^2\varphi_{R}^{2m-2}dxdy\right)^{1-\lambda}.
 \end{align*}
Therefore
\begin{align*}
\int_{ B_{R}\times B_{R^{1+s}}}\Big(1+\|\mathbf{x}\|^{2(s+1)}\Big)^{\frac{\alpha}{2(s+1)}} v^2dxdy \leq  CR^{N_{s}-\frac{2(\theta+1)p}{p\theta-1}-\frac{(p+1)\alpha}{p\theta-1} - \frac{2(2+\alpha)(2 - p)}{\theta+p-2}},
 \end{align*}
so we are done. \qed

\medskip

\section{Proofs of Theorem  \ref{main3} and Corollary \ref{main2}.}
\setcounter{equation}{0}

\medskip

Adopting the similar approach as in  Lemma 3.1 in  \cite{hhy}, we  establish the following lemma which plays an important role in proving
 Theorems \ref{main3}.

\begin{lem}
\label{newl}
Let $\alpha\geq0.$ Assume that $(u,v)$ is a stable solution of \eqref{1.1} such that $u$ is bounded. Then for any $z> \frac{p+1}{2}$ verifying
$L(z) < 0$, there exists $C <\infty$ such that
\begin{equation}
\label{3.1}
\int_{ B_{R}\times B_{R^{1+s}}}\Big(1+\|\mathbf{x}\|^{2(s+1)}\Big)^{\frac{\alpha}{2(s+1)}} u^{\theta}v^{z-1} dxdy \leq\frac{C}{R^2}\int_{ B_{2R}\times B_{(2R)^{1+s}}}v^{z} dxdy.
\end{equation}
where \begin{align}\label{L}
 L(z):=z^4-\frac{16p\theta(p+1)}{\theta+1}z^2+\frac{16p\theta(p+1)(p+\theta+2)}{(\theta+1)^2}z-\frac{16p\theta(p+1)^2}{(\theta+1)^2}.
\end{align}
\end{lem}
\noindent {\bf Proof}. Let $(u,v)$ be a stable solution of
\eqref{1.1}. Let $ \phi \in C_c^\infty\left(\mathbb{R}^N=\mathbb{R}^{N_1}\times\mathbb{R}^{N_2}, [0,1]\right),$ and $\varphi =
u^{\frac{q+1}{2}}\phi$ with $q > 0$. Integrating by parts, we get

\begin{align}\label{3.2}
 \begin{split}
\int_{\mathbb{R}^N}|\nabla_{s}\varphi|^2dxdy&=\frac{(q+1)^2}{4q}\int_{\mathbb{R}^N}\Big(1+\|\mathbf{x}\|^{2(s+1)}\Big)^{\frac{\alpha}{2(s+1)}}u^qv^p\phi^2dxdy+\int_{\mathbb{R}^N}u^{q+1}|\nabla_{s}\phi|^2dxdy
\\
&+\frac{1-q}{4q}\int_{\mathbb{R}^N}u^{q+1}\Delta_{s}(\phi^2)dxdy.
\end{split}
\end{align}
 Take $\varphi$ into the stability inequality \eqref{1.3} and using \eqref{3.2}, we obtain
\begin{align*}
&\sqrt{p\theta}\int_{\mathbb{R}^N}\Big(1+\|\mathbf{x}\|^{2(s+1)}\Big)^{\frac{\alpha}{2(s+1)}}u^{\frac{\theta-1}{2}}v^{\frac{p-1}{2}}u^{q+1}\phi^2dxdy\\
&\leq \int_{\mathbb{R}^N}|\nabla_{s}\varphi|^2dxdy\;\leq\; \frac{(q+1)^2}{4q}\int_{\mathbb{R}^N}\Big(1+\|\mathbf{x}\|^{2(s+1)}\Big)^{\frac{\alpha}{2(s+1)}}u^qv^p\phi^2dxdy+C\int_{\mathbb{R}^N}u^{q+1}\Big[|\nabla_{s}\phi|^2+\Delta_{s}(\phi^2)\Big]dxdy,
\end{align*}
so we get
\begin{align*}
& \;a_1\int_{\mathbb{R}^N}\Big(1+\|\mathbf{x}\|^{2(s+1)}\Big)^{\frac{\alpha}{2(s+1)}}u^{\frac{\theta-1}{2}}v^{\frac{p-1}{2}}u^{q+1}\phi^2dxdy\\
   \leq & \;
\int_{\mathbb{R}^N}\Big(1+\|\mathbf{x}\|^{2(s+1)}\Big)^{\frac{\alpha}{2(s+1)}}u^qv^p\phi^2dxdy +C\int_{\mathbb{R}^N}u^{q+1}\Big[|\nabla_{s}\phi|^2+\Delta_{s}(\phi^2)\Big]dxdy,
\end{align*}
where $a_1=\frac{4q\sqrt{p\theta}}{(q+1)^2}$. Choose now  $\phi(x,y)=\psi(\frac{x}{R},\frac{y}{R^{1+s}}),$  where $\psi $ a cut-off function in $ C_c^\infty\left(\mathbb{R}^N=\mathbb{R}^{N_1}\times\mathbb{R}^{N_2}, [0,1]\right),$ such that $$\psi=1 \quad \mbox{on} \quad B_{1}\times B_{1},\quad \mbox{and} \quad  \psi=0 \quad \mbox{outside } \quad B_{2}\times B_{2^{1+s}}.$$

A simple calculation,  implies that

$$|\nabla_{s}\phi|\leq \frac{C}{R}\quad \mbox{and} \quad  |\D_{s}(\phi^2)|\leq \frac{C}{R^{2}}.$$

Hence,
\begin{align*}
 \begin{split}
&I_1:=\int_{\mathbb{R}^N}\Big(1+\|\mathbf{x}\|^{2(s+1)}\Big)^{\frac{\alpha}{2(s+1)}}u^{\frac{\theta-1}{2}}v^{\frac{p-1}{2}}u^{q+1}\phi^2dxdy\\
&\leq
\frac{1}{a_1}
\int_{\mathbb{R}^N}\Big(1+\|\mathbf{x}\|^{2(s+1)}\Big)^{\frac{\alpha}{2(s+1)}}u^qv^p\phi^2dxdy+\frac{C}{R^2}\int_{ B_{2R}\times B_{(2R)^{1+s}}}u^{q+1}dxdy.
\end{split}
\end{align*}

Furthermore,, using $v^{\frac{r+1}{2}}\phi,$  $r
> 0$ as test function in \eqref{1.3}. As above,  we get readily
\begin{align*}
 \begin{split}
&I_2 :=\int_{\mathbb{R}^N}\Big(1+\|\mathbf{x}\|^{2(s+1)}\Big)^{\frac{\alpha}{2(s+1)}}u^{\frac{\theta-1}{2}}v^{\frac{p-1}{2}}v^{r+1}\phi^2dxdy\\
&\leq
\frac{1}{a_2} \int_{\mathbb{R}^N}\Big(1+\|\mathbf{x}\|^{2(s+1)}\Big)^{\frac{\alpha}{2(s+1)}}u^\theta
v^r\phi^2dxdy+\frac{C}{R^2}\int_{ B_{2R}\times B_{(2R)^{1+s}}}v^{r+1}dxdy.
\end{split}
\end{align*}
 with $a_2=\frac{4r\sqrt{p\theta}}{(r+1)^2}$. Combining  the two  last inequalities,  we have then
 \begin{align}
 \label{3.5}
  \begin{split}
 & I_1+{a_2}^\frac{2(r+1)}{p+1} I_2\\
   &\leq \frac{1}{a_1}
\int_{\mathbb{R}^N}\Big(1+\|\mathbf{x}\|^{2(s+1)}\Big)^{\frac{\alpha}{2(s+1)}}u^qv^p\phi^2dxdy+{a_2}^\frac{2r+1-p}{p+1}\int_{\mathbb{R}^N}\Big(1+\|\mathbf{x}\|^{2(s+1)}\Big)^{\frac{\alpha}{2(s+1)}}u^\theta
v^r\phi^2dxdy\\
   \;&\;+\frac{C}{R^2}\int_{ B_{2R}\times B_{(2R)^{1+s}}}\left(u^{q+1} +
v^{r+1}\right)dxdy.
  \end{split}
 \end{align}
Fix
\begin{equation}\label{3.6}
 q=\frac{(\theta+1)r}{p+1}+\frac{\theta-p}{p+1}, \quad \mbox{ or equivalently } \quad
q+1=\frac{(\theta+1)(r+1)}{p+1}.
\end{equation}
Let  $r> \frac{p-1}{2}.$  Applying Young's inequality and  using \eqref{3.6},  the  first term on the right hand side of \eqref{3.5}, can be estimated as
\begin{align*}
& \frac{1}{a_1}\int_{\mathbb{R}^N}\Big(1+\|\mathbf{x}\|^{2(s+1)}\Big)^{\frac{\alpha}{2(s+1)}}u^qv^p\phi^2dxdy\\
= & \; \frac{1}{a_1}\int_{\mathbb{R}^N}\Big(1+\|\mathbf{x}\|^{2(s+1)}\Big)^{\frac{\alpha}{2(s+1)}}u^{\frac{\theta-1}{2}}v^{\frac{p-1}{2}}u^{\frac{(\theta+1)r}{p+1}+\frac{\theta+1}{p+1}\left(\frac{1-p}{2}\right)}v^{\frac{p+1}{2}}\phi^2 dxdy \\
= & \;\frac{1}{a_1}\int_{\mathbb{R}^N}\Big(1+\|\mathbf{x}\|^{2(s+1)}\Big)^{\frac{\alpha}{2(s+1)}}u^{\frac{\theta-1}{2}}v^{\frac{p-1}{2}}u^{(q+1)\frac{2r+1-p}{2(r+1)}}v^{\frac{p+1}{2}}\phi^2dxdy\\
\leq & \; \frac{2r+1-p}{2(r+1)}\int_{\mathbb{R}^N}\Big(1+\|\mathbf{x}\|^{2(s+1)}\Big)^{\frac{\alpha}{2(s+1)}}u^{\frac{\theta-1}{2}}v^{\frac{p-1}{2}}u^{q+1}\phi^2dxdy\\
 &+ \;
 \frac{p+1}{2(r+1)}a_1^{-\frac{2(r+1)}{p+1}}\int_{\mathbb{R}^N}\Big(1+\|\mathbf{x}\|^{2(s+1)}\Big)^{\frac{\alpha}{2(s+1)}}u^{\frac{\theta-1}{2}}v^{\frac{p-1}{2}}v^{r+1}\phi^2dxdy
\\
= & \; \frac{2r+1-p}{2(r+1)}I_1+\frac{p+1}{2(r+1)} a_1^{-\frac{2(r+1)}{p+1}} I_2,
\end{align*}
and similarly
\begin{align*}
  {a_2}^{\frac{2r+1-p}{p+1}}\int_{\mathbb{R}^N}\Big(1+\|\mathbf{x}\|^{2(s+1)}\Big)^{\frac{\alpha}{2(s+1)}} u^\theta v^r\phi^2dxdy  \leq \frac{p+1}{2(r+1)} I_1
  + \frac{2r+1-p}{2(r+1)}{a_2}^\frac{2(r+1)}{p+1}I_2.
   \end{align*}

Inserting the two above estimates in \eqref{3.5}, we arrive at

\begin{align*}
{a_2}^{\frac{2(r+1)}{p+1}}I_2\leq
\left[\frac{2r+1-p}{2(r+1)}{a_2}^{\frac{2(r+1)}{p+1}}+\frac{p+1}{2(r+1)}{a_1}^{\frac{-2(r+1)}{p+1}}\right]I_2 +
\frac{C}{R^2}\int_{ B_{2R}\times B_{(2R)^{1+s}}}\left(u^{q+1}
+ v^{r+1}\right)dxdy.
\end{align*}
Combining \eqref{3.6} and \eqref{estS}, one obtains
$$u^{q+1} \leq Cv^{r+1}\quad \mbox{and} \quad
u^{\frac{\theta-1}{2}}v^{\frac{p-1}{2}}v^{r+1}\geq u^\theta v^r.$$
We get then
\begin{align*}
\frac{p+1}{2(r+1)}\left[(a_1a_2)^{\frac{2(r+1)}{p+1}}-1\right]
\int_{\mathbb{R}^N}u^\theta v^r\phi^2dxdy\leq CR^{-2}a_1^{\frac{2(r+1)}{p+1}}\int_{ B_{2R}\times B_{(2R)^{1+s}}}
v^{r+1}dxdy.
\end{align*}

Thus, if $a_1a_2 > 1$, there holds

\begin{align*}
\int_{ B_{R}\times B_{(R)^{1+s}}}\Big(1+\|\mathbf{x}\|^{2(s+1)}\Big)^{\frac{\alpha}{2(s+1)}}u^\theta v^r dxdy  \leq \int_{\mathbb{R}^N}u^\theta v^r\phi^2dxdy \leq \frac{C}{R^2}\int_{ B_{2R}\times B_{(2R)^{1+s}}}
v^{r+1}dxdy.
\end{align*}

Denote $z=r+1$, we conclude that if $a_1a_2
> 1$ and $z > \frac{p+1}{2}$,
\begin{align*}
 \int_{ B_{R}\times B_{(R)^{1+s}}}\Big(1+\|\mathbf{x}\|^{2(s+1)}\Big)^{\frac{\alpha}{2(s+1)}} u^\theta v^{z-1}dxdy  \leq \frac{C}{R^2}\int_{ B_{2R}\times B_{(2R)^{1+s}}} v^{z}dxdy.
\end{align*}
Furthermore, we can check that $a_{1}a_2>1$ is equivalent to $L(z)<0$, the proof is completed. \qed

\medskip

 Performing the change of variables $x=\frac{\theta+1}{p\theta-1}z$ in \eqref{newH}, a direct computation shows that
 $$H(x)=\left(\frac{\theta+1}{p\theta-1}\right)^4(z),$$
 where $L$ is given by \eqref{L}. Hence  $ H(x)<0$ if and only if $L(z)<0$. In addition, using  Lemma 6 in  \cite{Hfh}, we have

\begin{rem}\label{rnew2}
  \begin{itemize}
\item Let $1< p \leq \theta$, then $L(2) < 0$ and $L$ has a unique root $z_0$ in $(2, \infty)$ and $2t_0^+ < z_0$.
\item If $p > \frac{4}{3}$, then $L(p) < 0$ and $z_0$ is the unique root of $L$ in $(p, \infty),$ hence $x_0=\frac{\theta+1}{p\theta-1}z_0$.
\item  Therefore,  from Remark 3  in \cite{cow}, we find that$$
x_0 > 2t_0^+\frac{\theta+1}{p\theta-1}>4,\quad \forall\,\theta\geq p>1.$$
\end{itemize}
\end{rem}
\smallskip

We need to recall the following  properties of $t_{0}^+$ and $t_{0}^-$ before we complete the proof of  Theorem \ref{main3}.

\begin{rem}\label{rnew}
\begin{itemize}
\item It is known that for $1<p\leq \theta,$ there hold
$t_{0}^-<1<t_{0}^+,$ $t_{0}^-$ is decreasing and $t_{0}^+$ is increasing in $\varpi:=\frac{p\theta(p+1)}{\theta+1}.$ Moreover,
$\lim_{\varpi\rightarrow \infty}t_0^-= \frac{1}{2}$ and $\lim_{\varpi\rightarrow \infty}t_0^+= 1.$
\item Obviously $2t_0^-<p$ if $p>\frac{4}{3}$. Indeed, if $p>\frac{4}{3}$ then $\theta\geq p>\frac{4}{3}$ and $\varpi>\frac{16}{9}.$
Since $f(\varpi):=\sqrt{\varpi}-\sqrt{\varpi-\sqrt{\varpi}} $ is decreasing in $\varpi,$ there holds $2t_0^-=2f(\varpi)<2f(\frac{16}{9})=\frac{4}{3} < p.$
\end{itemize}
\end{rem}
\medskip\noindent
\subsection{\bf End of the proof of Theorem \ref{main3}.}
 \medskip

 In this subsection, we use $L^{2}$-estimates for Grushin operator,  and we apply the bootstrap iteration  as in \cite{DP,Hfh,cow}. For the completeness, we present the details.

\smallskip
 Let $ \eta \in C_c^\infty\left(\mathbb{R}^N=\mathbb{R}^{N_1}\times\mathbb{R}^{N_2}, [o,1]\right),$ be  a cut-off function  such that

 \begin{align}
\label{test}
\eta=1 \quad \mbox{on} \quad B_{1}\times B_{1},\quad \mbox{and} \quad  \eta=0 \quad \mbox{outside } \quad B_{2}\times B_{2^{1+s}}.
\end{align}

\smallskip

We divide the proof in three parts.

\medskip
{\bf Step $1.$} Denote by $\lambda_{s}=\frac{N_{s}}{N_{s}-2}.$  We claim that  for any smooth non-negative function $w,$ there exists  a positive constant $C>0$ such that

 \begin{align}
\label{test2}
\begin{split}
&\left(\int_{ B_{R}\times B_{R^{1+s}}} w^{2\lambda_{s}} dxdy\right)^{\frac{1}{\lambda_{s}}}\\
&\leq
CR^{N_{s}\big(\frac{1}{\lambda_{s}}-1\big)+2}\int_{ B_{2R}\times B_{(2R)^{1+s}}}|\nabla_{s} w|^{2}dxdy+CR^{N_{s}\big(\frac{1}{\lambda_{s}}-1\big)}\int_{ B_{2R}\times B_{(2R)^{1+s}}}w^{2}dxdy.
\end{split}
\end{align}

In Fact, By using Sobolev inequality  \cite{xxy} and integration by parts, imply that
\begin{align*}
\left(\int_{ B_{1}\times B_{1}} w^{2\lambda_{s}} dxdy\right)^{\frac{1}{2\lambda_{s}}}&\leq \left(\int_{ B_{2}\times B_{2^{1+s}}} (w\eta)^{2\lambda_{s}} dxdy\right)^{\frac{1}{2\lambda_{s}}}\\
&\leq C\left(\int_{ B_{2}\times B_{2^{1+s}}} |\nabla_{s} (w\eta)|^{2} dxdy\right)^{\frac{1}{2}}\\
&\leq C\left[\int_{ B_{2}\times B_{2^{1+s}}} \left(|\nabla_{s} w|^{2}\eta^{2}+w^{2}|\nabla_{s} \eta|^{2}-\frac{w^{2}}{2}\D_{s}(\eta)  \right)dxdy\right]^{\frac{1}{2}},
\end{align*}
So, we get
\begin{align*}
\left(\int_{ B_{1}\times B_{1}} w^{2\lambda_{s}} dxdy\right)^{\frac{1}{\lambda_{s}}}
\leq C\int_{ B_{2}\times B_{2^{1+s}}} \left(|\nabla_{s} w|^{2}+w^{2}  \right)dxdy.
\end{align*}
By scaling argument, we obtain  readily the estimate \eqref{test2}.

\medskip
{\bf Step $2.$} Let $(u,v)$ be a  stable solution of \eqref{1.1}, with $1< p \leq \theta$. Then for any  $\lambda_{s}=\frac{N_{s}}{N_{s}-2}$ and $2t_0^-<z_0$, we claim that there exists a positive constant $C>0$ such that
\begin{equation}
\label{lnew2}
\left(\int_{ B_{R}\times B_{R^{1+s}}} v^{z_0\lambda_{s}} dxdy\right)^{\frac{1}{\lambda_{s}}}\leq
CR^{N_{s}\big(\frac{1}{\lambda_{s}}-1\big)}\int_{ B_{2R}\times B_{(2R)^{1+s}}}v^{z_0}dxdy.
\end{equation}

To prove this, for $2t_0^-<z_0,$ in what follows, we choose
$$w=v^{\frac{z_0}{2}}.$$

Let us put $\eta_{R}(x,y)=\eta(\frac{x}{R},\frac{y}{R^{1+s}}),$ where $\eta$ is given in \eqref{test}. A simple calculation, we obtain readily
 \begin{align}
\label{stbk}
\int_{ B_{R}\times B_{R^{1+s}}} |\nabla_{s} w|^{2} dxdy \leq C \int_{ B_{2R}\times B_{(2R)^{1+s}}}v^{z_{o}-2}|\nabla_{s} v|^{2}\eta_{R}^{2} dxdy.
\end{align}
 Multiplying $-\D v= \Big(1+\|\mathbf{x}\|^{2}\Big)^{\frac{\alpha}{2}}u^{\theta}$ by $v^{z_{o}-1} \eta_{R}^2$ and  integrating by parts, we derive
\begin{align}\label{3.9}
\begin{split}
 (z_{o}-1) \int_{ B_{2R}\times B_{(2R)^{1+s}}}v^{z_{o}-2}|\nabla_{s} v|^{2}\eta_{R}^2dxdy
&\;= \int_{ B_{2R}\times B_{(2R)^{1+s}}}\Big(1+\|\mathbf{x}\|^{2(s+1)}\Big)^{\frac{\alpha}{2(s+1)}}v^{z_{o}-1} u^{\theta}\eta_{R}^2 dxdy\\
&\;-2 \int_{ B_{2R}\times B_{(2R)^{1+s}}}\eta_{R} v^{z_{o}-1}\nabla_{s}v \cdot \nabla_{s} \eta_{R}  dxdy.
\end{split}
\end{align}
By Young's inequality,
\begin{align*}
 &\; 2 \int_{ B_{2R}\times B_{(2R)^{1+s}}} v^{z_{0}-1}|\nabla_{s} v||\nabla_{s} \eta_{R}|\eta_{R} dxdy \\
&\; \leq \frac{z_{0}-1}{2} \int_{ B_{2R}\times B_{(2R)^{1+s}}}v^{z_{0}-2} |\nabla_{s} v|^{2}\eta_{R}^2dxdy +
  C \int_{ B_{2R}\times B_{(2R)^{1+s}}}v^{z_{0}}|\nabla_{s} \eta_{R}| ^2dxdy.
\end{align*}
Substituting this in \eqref{3.9}, and by virtue of estimate  \eqref{stbk},  we arrive at
\begin{align*}
&\;\int_{ B_{R}\times B_{R^{1+s}}} |\nabla_{s} w|^{2} dxdy\\
&\;\leq \int_{ B_{2R}\times B_{(2R)^{1+s}}}v^{z_{o}-2}|\nabla_{s} v|^{2}\eta_{R}^2dxdy \\
&\;\leq C\int_{ B_{2R}\times B_{(2R)^{1+s}}}\Big(1+\|\mathbf{x}\|^{2(s+1)}\Big)^{\frac{\alpha}{2(s+1)}}v^{z_{o}-1} u^{\theta}\eta_{R}^2 dxdy+\frac{C}{R^{2}}\int_{ B_{2R}\times B_{(2R)^{1+s}}}v^{z_{0}}dxdy.
\end{align*}
In view of estimate \eqref{test2}, and using Lemmas \ref{newl}, we obtain easily
the estimate \eqref{lnew2}.

\medskip
{\bf Step $3.$} Let $(u,v)$ be a  stable solution of \eqref{1.1}, with $1< p \leq \theta$  and $q \in (2t_0^-, z_0).$ Then for any  $\lambda_{s}=\frac{N_{s}}{N_{s}-2},$ $q<z_m\lambda_{s},$ we claim that there exists a positive constant $C>0$ such that
 \begin{align}
 \label{lnew2j}
 \left(\int_{ B_{R}\times B_{R^{1+s}}} v^{z_m\lambda_{s}} dxdy\right)^{\frac{1}{z_m\lambda_{s}}}\leq
CR^{N_{s}\big(\frac{1}{z_m\lambda_{s}}-\frac{1}{q}\big)}\left(\int_{ B_{R_{m}}\times B_{(R_{m})^{1+s}}}v^{q}dxdy\right)^{\frac{1}{q}}.
\end{align}

To prove this, by Remark \ref{rnew}, we known that $2t_0^-<p.$ Fix a real positive number $q$ satisfying
 $$2t_0^-<q<p.$$
 Let $m$ be the nonnegative integer such that $q\lambda_{s}^{m-1}<z_{0}<q\lambda_{s}^{m}.$ Put
 $$z_{1}=qk,  \quad z_{2}=qk\lambda_{s},....,z_{m}=qk\lambda_{s}^{m-1},$$
 verifying
 $$2t_0^-<z_{1}<z_{2}<,....,<z_{m}<z_{0}.$$
 Here the constant $k\in[1,\lambda_{s}],$  is chosen such that $z_{m}$ is arbitrarily close to  $z_{0}.$
 Set $ R_n = 2^n R.$ by \eqref{lnew2}, and an induction argument, we deduce then
 \begin{align}
\label{lnnnew2}
\begin{split}
\left(\int_{ B_{R}\times B_{R^{1+s}}} v^{z_m\lambda_{s}} dxdy\right)^{\frac{1}{z_m\lambda_{s}}}&\leq
CR^{N_{s}\big(\frac{1}{z_m\lambda_{s}}-\frac{1}{z_m}\big)}\left(\int_{ B_{1}\times B_{(1)^{1+s}}}v^{z_m}dxdy\right)^{\frac{1}{z_m}}\\
&\;=CR^{N_{s}\big(\frac{1}{z_m\lambda_{s}}-\frac{1}{z_m}\big)}\left(\int_{ B_{1}\times B_{(1)^{1+s}}}v^{z_{m-1}\lambda_{s}}dxdy\right)^{\frac{1}{z_{m-1}\lambda_{s}}}\\
&\;\leq CR^{N_{s}\big(\frac{1}{z_m\lambda_{s}}-\frac{1}{z_1}\big)}\left(\int_{ B_{R_{m}}\times B_{(R_{m})^{1+s}}}v^{z_{1}}dxdy\right)^{\frac{1}{z_{1}}}\\
&\;\leq CR^{N_{s}\big(\frac{1}{z_m\lambda_{s}}-\frac{1}{qk}\big)}\left(\int_{ B_{R_{m}}\times B_{(R_{m})^{1+s}}}v^{qk}dxdy\right)^{\frac{1}{qk}}
\end{split}
\end{align}
 Furthermore, by H\"older's inequality, there holds
  \begin{align}
\label{lnnnew23}
\begin{split}
&\;\left(\int_{ B_{R_{m}}\times B_{(R_{m})^{1+s}}}v^{qk}dxdy\right)^{\frac{1}{qk}}\\
&\;\leq
\left[\left(\int_{ B_{R_{m}}\times B_{(R_{m})^{1+s}}}v^{q\lambda_{s}}dxdy\right)^{\frac{k}{\lambda_{s}}}\left(\int_{ B_{R_{m}}\times B_{(R_{m})^{1+s}}}dxdy\right)^{1-\frac{k}{\lambda_{s}}}\right]^{\frac{1}{qk}}\\
&\;\leq C\left[\left(\int_{ B_{R_{m}}\times B_{(R_{m})^{1+s}}}v^{q\lambda_{s}}dxdy\right)^{\frac{k}{\lambda_{s}}}CR^{N_{s}\big(1-\frac{k}{\lambda_{s}}\big)}\right]^{\frac{1}{qk}}\\
&\;\leq CR^{N_{s}\big(\frac{1}{kq}-\frac{1}{q\lambda_{s}}\big)}\left(\int_{ B_{R_{m}}\times B_{(R_{m})^{1+s}}}v^{q\lambda_{s}}dxdy\right)^{\frac{1}{q\lambda_{s}}}\\
&\;\leq CR^{N_{s}\big(\frac{1}{kq}-\frac{1}{q\lambda_{s}}\big)}R^{N_{s}\big(\frac{1}{q\lambda_{s}}-\frac{1}{q}\big)}\left(\int_{ B_{R_{m}}\times B_{(R_{m})^{1+s}}}v^{q}dxdy\right)^{\frac{1}{q}}.
\end{split}
\end{align}
 Therefore, the claim follows by combining the last tow inequality.

\medskip\noindent
{\bf Proof of Theorem \ref{main3} completed.} We are now in position to conclude.  Let $\alpha\geq0,$  and $(u, v)$ be a  stable solution of \eqref{1.1} with $\theta \geq p > 1.$   We split the proof into two cases:
$p >\frac{4}{3}$ and  $1< p \leq \min(\frac{4}{3}, \theta).$

\medskip
\textit{Case 1: $p >\frac{4}{3}$.} Let $p>q>0.$ From \eqref{VV}, we use H\"older's inequality, to obtain

\begin{align}
\label{3.13}
\begin{split}
\int_{ B_{R}\times B_{R^{1+s}}}v^qdxdy
&\leq \left(\int_{ B_{R}\times B_{R^{1+s}}}\Big(1+\|\mathbf{x}\|^{2(s+1)}\Big)^{\frac{\alpha}{2(s+1)}} v^pdxdy\right)^{\frac{q}{p}} \\
&\; \;\;\;\times\left(\int_{ B_{R}\times B_{R^{1+s}}}\Big(1+\|\mathbf{x}\|^{2(s+1)}\Big)^{-\frac{\alpha p}{2(s+1)(p-q)}}dxdy\right)^{\frac{p-q}{p}}\\
&\;\leq C R^{\left[N_{s}-\frac{2(\theta+1)p}{p\theta-1}-\frac{(p+1)\alpha}{p\theta-1}\right]\frac{q}{p}
+\left(N_{s}-\frac{\alpha q}{p-q}\right)\frac{p-q}{p}}
         = CR^{N_{s}-\frac{(2+\alpha)(\theta+1)}{p\theta-1}q}.
\end{split}
\end{align}

Combining \eqref{3.13} and  \eqref{lnew2j}, we deduce that

\begin{align}
 \label{lmbv}
 \left(\int_{ B_{R}\times B_{R^{1+s}}} v^{z_m\lambda_{s}} dxdy\right)^{\frac{1}{z_m\lambda_{s}}}\leq
CR^{\frac{N_{s}}{z_m\lambda_{s}}-\frac{(2+\alpha)(\theta+1)}{p\theta-1}}.
\end{align}
 Recall that $\lambda_{s}=\frac{N_{s}}{N_{s}-2}.$ Since
$$N<2+\left(\frac{(2+\alpha)(\theta+1)}{p\theta-1}\right)z_0,$$

 we chose $k\in[1,\lambda_{s}],$ such that $z_{m}$ is close to  $z_{0}.$ Then, it implies from \eqref{lmbv} that $\|v\|_{L^{z_m\lambda_{s}}(\R^N)} = 0$ as $R\rightarrow\infty,$ i.e., $v \equiv 0$ in $\R^N$ This is a contraction.
Therefore, we get the desired result: the equation \eqref{1.1} has no 
stable solution if $N_{s}<2+(2+\alpha)x_0$ where $x_0=\frac{\theta+1}{p\theta-1}z_0$.

\smallskip
Finally, it follows from   Remark \ref{rnew2} that if $N_{s}\leq 10+4\alpha$, \eqref{1.1} has no  stable solution
for any $\theta \geq p > \frac{4}{3}.$

\medskip
\textit{Case 2: $1<p\leq \frac{4}{3}$ and $u$ is bounded.} Let now $2>q>0,$ using  Lemma \ref{lemnew}, we obtain
\begin{align*}
&\int_{ B_{R}\times B_{R^{1+s}}}v^qdxdy\\
&\;\leq \left(\int_{ B_{R}\times B_{R^{1+s}}}\Big(1+\|\mathbf{x}\|^{2(s+1)}\Big)^{\frac{\alpha}{2(s+1)}} v^2dxdy\right)^{\frac{q}{2}} \left(\int_{ B_{R}\times B_{R^{1+s}}}\Big(1+\|\mathbf{x}\|^{2(s+1)}\Big)^{-\frac{\alpha q}{2(s+1)(2-q)}}dxdy\right)^{\frac{2-q}{2}}\\
& \leq C R^{\left[N_{s}-\frac{2(\theta+1)p}{p\theta-1}-\frac{(p+1)\alpha}{p\theta-1} - \frac{2(2+\alpha)(2 - p)}{\theta+p-2}\right]\frac{q}{2}+\left(N_{s}-\frac{\alpha q}{2-q}\right)\frac{2-q}{2}}\\
& = CR^{N_{s}-\left[\frac{(\theta+1)p}{p\theta-1} + \frac{(2+\alpha)(2 - p)}{\theta+p-2} +\frac{p(\theta+1)\alpha}{2(p\theta-1)}\right] q},
      \end{align*}
      Substituting this in  \eqref{lnew2j},
\begin{align*}
 \left(\int_{ B_{R}\times B_{R^{1+s}}} v^{z_m\lambda_{s}} dxdy\right)^{\frac{1}{z_m\lambda_{s}}}\leq
CR^{\frac{N_{s}}{z_m\lambda_{s}}-\frac{(\theta+1)p}{p\theta-1} + \frac{(2+\alpha)(2 - p)}{\theta+p-2} +\frac{p(\theta+1)\alpha}{2(p\theta-1)}}.
\end{align*}
Arguing as Case 1, we get the desired result \qed

\medskip\noindent
\subsection{\bf End of the proof of Corollary \ref{main2}.}

 \medskip
If $p = \theta>1,$  then, from Lemma \ref{Soup}, we get that $v = u$ and the weighted Lane-
Emden system \eqref{1.1} is reduced to be the weighted Lane-Emden equation \eqref{D}. From Remark \ref{rnew2} and arguing as above, we get
\begin{align*}
  L(z) = z^4 -16p^2z^2 +
32p^2z-16p^2=(z^2+4p(z-1))(z-2t_0^-)(z-2t_0^+) \quad \mbox{with } \quad t_0^{\pm}= p\pm\sqrt{p^2-p}.
\end{align*}

Adopting the similar argument as in the proof  of Corollary 1.2  in  \cite{Hfh}, we  obtain $2t_0^+$ is the largest root of $L$ as $t_0^+ > p > 1$. Therefore
$$x_0 = \frac{2p+2\sqrt{p^2-p}}{p-1}>4\quad \mbox{for all } \, p>1.$$  The result follows immediately by applying Theorem \ref{main3}. \qed

\bigskip\noindent
{\bf Acknowledgments}. We would like to thank Professor Dong Ye for suggesting us this problem and for many helpful comments.

\end{document}